\documentclass[11pt]{amsart}
\usepackage{amssymb,amsmath,txfonts,mathrsfs,mathtools,titletoc, tikz, stmaryrd}
\usepackage[alphabetic]{amsrefs}
\usepackage{float}
\usepackage[colorlinks, linkcolor=blue,anchorcolor=Periwinkle,
    citecolor=red,urlcolor=Emerald]{hyperref}
\usepackage{xcolor}
\usepackage{enumitem}
\usepackage{geometry,array} 
\usepackage{graphicx}
\usepackage{subfigure}
\usepackage{bookmark}
\usepackage{tikz}\usetikzlibrary{matrix}
\usepackage{url}
\usepackage{braids}

\usetikzlibrary{decorations.markings}
\tikzset{->-/.style={decoration={  markings,  mark=at position #1 with
    {\arrow{>}}},postaction={decorate}}}
\tikzset{-<-/.style={decoration={  markings,  mark=at position #1 with
    {\arrow{<}}},postaction={decorate}}}
\usepackage[all]{xy}
\usetikzlibrary{arrows}
\usetikzlibrary{decorations.pathreplacing,decorations.pathmorphing,shadings,fadings,calc}
\usepackage{extarrows}
\usepackage{tikz-3dplot}
\usepackage{appendix}

\usetikzlibrary{matrix}
\usepackage[all]{xy}
\usetikzlibrary{arrows,shapes}
\usetikzlibrary{calc}
\tikzset{->-/.style={decoration={  markings,  mark=at position #1 with
    {\arrow{>}}},postaction={decorate}}}
\tikzset{-<-/.style={decoration={  markings,  mark=at position #1 with
    {\arrow{<}}},postaction={decorate}}}
\usepackage[]{hyperref}
\numberwithin{equation}{section}

\usetikzlibrary{decorations.markings}
\tikzset{->-/.style={decoration={  markings,  mark=at position #1 with
    {\arrow{>}}},postaction={decorate}}}
\tikzset{-<-/.style={decoration={  markings,  mark=at position #1 with
    {\arrow{<}}},postaction={decorate}}}

\usepackage{lscape}
\usepackage{CJK}

\newtheorem{theorem}{Theorem}[section]

\newtheorem{claim}[theorem]{Claim}

\newtheorem{definition}[theorem]{Definition}

\numberwithin{equation}{section}

\def\QEDopen{{\setlength{\fboxsep}{0pt}\setlength{\fboxrule}{0.2pt}\fbox{\rule[0pt]{0pt}{1.3ex}\rule[0pt]{1.3ex}{0pt}}}}
\def\QED{\hfill\QEDopen} 
\def\proof{\noindent{\it Proof}: }

\newcommand{\R}{\mathbb{R}}

\begin{document}
\def\nn{node{$\bullet$}}
\def\ww{node{$\circ$}}

\title{The Distance Between the Perturbation of a Convex Function and its $\Gamma$-regularization}
\author[1]{Zichang Liu}
\address{Department of Mathematical Sciences\\Tsinghua University, Beijing\\P. R. China, 100084}
\email{liu-zc19@mails.tsinghua.edu.cn}
\date{\today}

\begin{abstract}
In the study of a non-convex minimization problem by Lachand-Robert and Peletier, they found that the difference between the compactly supported perturbation $u+\epsilon h$ of a strictly convex function $u$, and the $\Gamma$-regularization of $u+\epsilon h$, is at most $o(\epsilon)$. Here we find that this result is optimal, albeit they expected a much stronger estimate.
\end{abstract}
\thanks{}

\maketitle

\section{Introduction}
In \cite{lrp}, Lachand-Robert and Peletier studied the minimizer of the functional
\begin{equation}
    \int_\Omega f(\nabla u) dx, \label{1}
\end{equation}
where $\Omega \subset \R^2$ is a bounded open set, $f:\R^2\rightarrow \R$ is a non-negative smooth function, and the minimum is taken over
$$\mathcal{C}:=\{u: \Omega\rightarrow [0,1];\ u\ \text{convex}\}.$$
Unlike general variational problems, this problem cannot be studied using the usual methods of the calculus of variations, since $u$ stays in the class of convex functions. Moreover, even if $f$ is convex, there is no regularizing effect as the usual variational problems.\\
In particular, \cite{lrp} considered the case where the Hessian matrix $d^2 f$ has at least one negative eigenvalue at each point. In this case $f$ is said to be \textit{nowhere convex}. An insightful observation suggests that if $u$ is $C^2$ in an open subset of $\Omega$ with positive det$d^2 u$ in this set, then $u$ is not the desired minimizer. This observation suggests that a minimizer $u$ should satisfy ``det$d^2 u=0$" in some sense. Unfortunately, one can almost never expect the minimizer is of class $C^2$, since there is no regularizing effect in this sort of problems.\\
Actually, under some further technical assumptions on the function $f$, the authors proved the following theorem:
\begin{theorem}
    \label{thm1}
    Let $u$ be a minimizer of \eqref{1} in $\mathcal{C}$, and $\Omega_1$ be an open convex subset of $\Omega$. If $f$ is nowhere convex on the convex hull of $\nabla u(\Omega_1)$, then $u$ is not strictly convex on $\Omega_1$.
\end{theorem}
Theorem \ref{thm1} suggests that there is a dense collection of line segments in the graph of $u$; one may consider cylindrical or conical surfaces as examples.\\
The study of minimizing \eqref{1} arose from the problem of the body of minimal resistance, which was introduced by Sir Isaac Newton in \textit{Principia Mathematica}. The latter problem can be stated as seeking for minimizers of \eqref{1} with
$$f(x)=\frac{1}{1+|\nabla x|^2},$$
in the set 
$$\mathcal{C}_M:=\{u: \Omega\rightarrow [0,M];\ u\ \text{convex}\}$$
with $M>0$. When $\Omega$ is a disk, Newton computed a minimizer $u_{rad}$ of \eqref{1} in the subset of $\mathcal{C}_M$ consisting of all radially symmetric functions. However, it was shown in \cite{bfk} that the radially symmetric minimizer $u_{rad}$ does not actually minimize \eqref{1} in $\mathcal{C}_M$. In particular, the minimizer of \eqref{1} is not unique. This is due to the fact that the Hessian matrix of $\frac{1}{1+|x|^2}$ has a negative eigenvalue at every point; there is a ``non-radial" direction at $u_{rad}$ in which \eqref{1} has a negative second variation.\\
Unfortunately, the argument given in \cite{bfk} provides no information on the shape of
the minimizers other than the lack of radial symmetry. However, we can prove their result in an alternative way using Theorem \ref{thm1}: one simply notice that the radially symmetric minimizer $u_{rad}$ is strictly convex on some open subset of $\Omega$.\\
In the proof of Theorem \ref{thm1}, the authors used the following technical result (Lemma 2 in \cite{lrp}): if $u$ is strictly convex in $\Omega$ and $h$ is a Lipschitz function compactly supported in $\Omega$, then the $\Gamma$-regularization of $u_\epsilon=u+\epsilon h$, denoted by $\Tilde{u_\epsilon}$, satisfies the following estimate:
$$|u_\epsilon-\Tilde{u_\epsilon}|_{C^0(\Omega)}=o(\epsilon),\ \epsilon\rightarrow 0.$$
We first give the formal definition of ``$\Gamma$-regularization":
\begin{definition}
Assume $\Omega\subset \R^n$ is a bounded open set, and $u:\Omega\rightarrow \R$ is a bounded function. $u^{**}:\R^n\rightarrow \R$, the \textbf{$\Gamma$-regularization}  of $u$, is defined to be the supremum of all affine functions which are less than or equal to $u$ on $\Omega$; inside $\Omega$ $u^{**}$ equals to the largest convex function less than or equal to $u$.
\end{definition}
In fact, the authors expected that the optimal convergence rate is $O(\epsilon^2)$, or even $o(\epsilon^2)$, and claimed that the better estimates could simplify the proof of Theorem \ref{thm1} considerably. However this is false: we will show by a simple example that the convergence rate $o(\epsilon)$ is already sharp.

\section{Main Result}
\begin{theorem}
    Let $\Omega=B_0(1)\subset \R^2$ (the unit disk). Given any function $\phi:[0,\delta]\rightarrow \R$, $\delta>0$, satisfying $\phi(0)=0$, $\phi(t)>0$ ($0<t\leq \delta$), $\phi(t)\searrow 0$ as $t\searrow 0$, there exists a smooth, strictly convex function $u:\bar{\Omega}\rightarrow \R$ and a Lipschitz function $h:\Omega\rightarrow \R$ compactly supported in $\Omega$, satisfying $|\Omega-\Omega_c|=0$, where $\Omega_c$ is the union of all open subsets of $\Omega$ on which $h$ is convex, s.t.
    $$\liminf_{t\searrow 0}\frac{u_t(0)-(u_t)^{**}(0)}{t\phi(t)}>0,$$
    where $u_t=u+th$.
\end{theorem}

\proof
We may assume $\phi$ is continuous and strictly increasing (by enlarging the values of $\phi$ suitably); denote by $\phi^{-1}$ the inverse function of $\phi$.\\
Let $\psi:[-1,1]\rightarrow \R$ be a smooth even function satisfying $\psi(0)=0$, $\psi(x)>0$ ($x\neq 0$), and $\psi(x)=O(x \phi^{-1}(x))$ as $x\searrow 0$, and let $v:[-1,1]\rightarrow \R$ satisfy
$$\begin{cases}
    v''(x)=\psi(x),\\
    v(0)=v'(0)=0,
\end{cases}$$
then $v$ is smooth, even, strictly convex, and clearly $v(x)=o(x \phi^{-1}(x))$ as $x\searrow 0$.\\
Define $g:[-1,1]\rightarrow \R$ as
$$g(x)=\begin{cases}
    \frac{1}{2}-|x|,\ |x|\leq \frac{1}{2},\\
    0,\ |x|\geq \frac{1}{2}.
\end{cases}$$
For any $t>0$, $v_t:=v+tg$ attains its minimum $m_t=v_t(\pm x_t)$ at $\pm x_t$, where $0\leq x_t\leq 1$. Thus, the graph of $(v_t)^{**}$ must contain the line segment joining $(-x_t,m_t)$ and $(x_t,m_t)$, i.e. $(v_t)^{**}(0)=m_t$.\\
Now, for $0<t\ll 1$,
$$\frac{v_t(0)-(v_t)^{**}(0)}{t\phi(t)}=\frac{v_t(0)-m_t}{t\phi(t)}\geq \frac{v_t(0)-v_t(\phi(t))}{t\phi(t)}=\frac{\frac{1}{2}t-\{v(\phi(t))+t(\frac{1}{2}-\phi(t))\}}{t\phi(t)}$$
$$=\frac{t\phi(t)-v(\phi(t))}{t\phi(t)}\rightarrow 1\ (t\searrow 0),$$
since $v(\phi(t))=o(t\phi(t))$.\\
Let $u(x_1,x_2)=v(\sqrt{x_1^2+x_2^2})$, then $u$ is smooth and strictly convex in $\bar{\Omega}$. Let $h(x_1,x_2)=g(\max\{|x_1|,|x_2|\})$, then $h$ is Lipschitz, compactly supported in $\Omega$, and satisfies our condition about convexity. Write $u_t=u+th$. Noting the fact that $(u_t)^{**}(x,0)$ is a convex function of $x$ satisfying $(u_t)^{**}(x,0)\leq u_t(x,0)=v_t(x)$, we have (for $|x|\leq 1$) $(u_t)^{**}(x,0)\leq (v_t)^{**}(x)$, and in particular $(u_t)^{**}(0,0)\leq (v_t)^{**}(0)$. Thus,
$$\frac{u_t(0)-(u_t)^{**}(0)}{t\phi(t)}\geq \frac{v_t(0)-(v_t)^{**}(0)}{t\phi(t)},$$
and the desired result follows.

\QED

\noindent{\textbf{Acknowledgement}}\\

The author wish to express his gratitude to Professor Guoyi Xu for bringing this problem to his attention. The author is supported by National Key R$\&$D Program of China 2022YFA1005400 and NFSC No.12031017.

\end{document}